\input amstex1
\input amssym.def
\documentstyle{amsppt1}
\input epsf
\magnification=1200

\newsymbol\varkappa 207B
\font\titlefont=cmr10 scaled \magstep3

\def\today{{\ifcase\month\or
 January\or February\or March\or April\or May\or June\or
 July\or August\or September\or October\or November\or December\fi
 \space\number\day, \number\year}}

\def\Re{{\,\text{\rm Re}\,}}    
\def\D{{\Bbb D}}       \def\T{{\Bbb T}}       \def\OD{{\overline\D}}
\def\C{{\Bbb C}}       
           
\def\calP{\Cal P}

\def\R{{\Bbb R}}

\def\l{\lambda}         \def\o{\omega}

\def\wbergman{L^2_a(\o)}       
\def\ltwo{L^2(\o)}     \def\hbergman{L^2_h(\o)}

\def\Hed{\hbox{\bf 1}}
\def\Rud{\hbox{\bf 2}}
\def\Shi{\hbox{\bf 3}}
\TagsOnRight

\topmatter
\title\titlefont
C\lowercase{urious properties of canonical divisors}
\\
\titlefont\lowercase{in weighted} B\lowercase{ergman spaces}
\endtitle
\author
Alexandru Aleman, H\aa kan Hedenmalm, Stefan Richter, Carl Sundberg
\endauthor
%\shortau{  }
%\abstract\nofrills{}
\address{\smc
A.~Aleman, Fachbereich Mathematik, Fernuniversit\"at Hagen, D-58084 Hagen, 
Germany
\newline\sl
Alexandru.Aleman\@FernUni-Hagen.de
\newline\smc
P.~J.~H.~Hedenmalm, Department of Mathematics, Lund University, Box 118,
S-221 00 Lund, Sweden
\newline\sl
E-mail: haakan\@maths.lth.se
\newline\smc
Stefan Richter and Carl Sundberg, Department of Mathematics, The
University of Tennessee, Knoxville, TN 37996, USA
\newline\sl
E-mail: richter\@novell.math.utk.edu
}
\date{
\vbox{\vskip 0.4 cm
\centerline{April, 1996}}
}

\endtopmatter

\document

\heading
1. Introduction 
\endheading
Let $\o$ be a nonnegative Borel measurable function in the open unit disk 
$\D$, such that
$$h(0)=\int_\D h(z)\o(z)dS(z),\qquad h\in L^\infty_h(\D),\tag1-1$$
where $dS$ stands for area measure in the complex plane, normalized so that
$\D$ has mass $1$, and $L^\infty_h(\D)$ is the space of bounded 
complex-valued harmonic functions on $\D$. We associate with $\o$ the Hilbert
space $L^2(\o)=L^2(\D,\o)$ of Borel measurable functions $f$ on $\D$ with
$$\|f\|_{L^2(\o)}=\left(\int_\D|f(z)|^2\o(z)dS(z)\right)^{1/2}<+\infty.$$
Consider the subspaces $\wbergman$ and $\hbergman$ of $\ltwo$ consisting of
those functions that can be altered on null sets so as to be analytic and
harmonic on $\D$, respectively. If two harmonic functions on $\D$ are given,
which coincide as elements of $\ltwo$, then the harmonic functions coincide.
This allows us to think of $\wbergman$ and $\hbergman$ as spaces of analytic
and harmonic functions, respectively. These linear subspaces of $\ltwo$ need
not be closed. For instance, if $\o$ is supported on a disk centered at the
origin of radius less than $1$, they certainly are not. However, if $\o$ is
locally bounded away from $0$ on $\D\setminus K$, for some compact subset $K$
of $\D$, then they are, in which case we refer to them as weighted analytic
and harmonic Bergman spaces, with weight $\o$. We sometimes drop the word
``analytic'', and talk about weighted Bergman spaces. When the weight is the
constant function $\o=1$, we talk about the analytic and harmonic Bergman
spaces, and write $L^2_a$ and $L^2_h$, dropping the weight from the notation.
In the sequel we shall make the assumption that $\o$ is such that $L^2_a(\o)$
and $L^2_h(\D)$ are closed subspaces of $L^2(\o)$.

Condition \thetag{1-1} implies that the decomposition of a harmonic polynomial
$f$ into a sum of an analytic polynomial $g$ and an antianalytic polynomial
$h$, where $h(0)=0$, is actually orthogonal: $g$ and $h$ are perpendicular in
$L^2_h(\o)$. If the harmonic polynomials are dense in $L^2_h(\o)$, this 
extends to a direct sum decomposition $L^2_h(\o)=L^2_a(\o)\oplus L^2_{\bar a,
0}(\o)$, where $L^2_{\bar a,0}(\o)$ stands for the closed subspace of $L^2_h
(\o)$ consisting of antianalytic functions that vanish at the origin.
In the sequel we shall assume that $\o$ is such that the harmonic polynomials
are dense in $L^2_h(\o)$.

\heading 2. Basic concepts 
\endheading 

The reproducing kernel function in any of the spaces $L^2_h(\o)$, $L^2_a(\o)$,
and $L^2_{\bar a,0}(\o)$, is obtained by taking an orthonormal basis $e_1,e_2,
\ldots$, and forming
$$k(z,\zeta)=\sum_{n=1}^\infty e_n(z)\overline{e}_n(\zeta),\qquad (z,\zeta)\in
\D\times\D.$$
To distinguish the kernels for the various spaces we use subscripts and 
superscripts: $k^\o_h(z,\zeta)$, $k^\o_a(z,\zeta)$, and $k^\o_{\bar a,0}(z,
\zeta)$, with obvious interpretations. It should be pointed out that the sum 
is independent of the particular orthonormal basis. One way to see this is
to use the reproducing property of the kernel functions. For instance, $k_h^\o
(z,\zeta)$ has the property that
$$f(z)=\int_\D f(\zeta)k^\o_h(z,\zeta)\o(\zeta)dS(\zeta),\qquad z\in\D,$$
for all $f\in L^2_h(\o)$. The direct decomposition of $L^2_h(\o)$ has its
counterpart for the kernel functions,
$$k_h^\o(z,\zeta)=k_a^\o(z,\zeta)+k_{\bar a,0}^\o(z,\zeta),\qquad (z,\zeta)\in
\D\times\D.$$
One checks that $k_{\bar a,0}^\o(z,\zeta)=\overline{k}_a^\o(z,\zeta)-1$, so
that the above identity reduces to
$$k_h^\o(z,\zeta)=2\Re k_a^\o(z,\zeta)-1,\qquad (z,\zeta)\in\D\times\D.
\tag1-2$$
We shall work with weights $\o$ such that the kernel $k_a^\o(z,\zeta)$ extends
continuously to ${\overline{\D}}^2\setminus\delta(\T)$, where $\delta(\T)=
\big\{(z,z):\,z\in\T\big\}$ is the diagonal on $\T^2$. If, for instance, 
$\o$ is $C^\infty$-smooth and bounded away from $0$ near $\T$, this is so. 

\heading 3. General properties of one-point extremal functions
\endheading

We assume the weight $\o$ is $C^\infty$ smooth near $\T$ in $\OD$, and bounded
away from $0$ there too. Then the kernel $k_a^\o(z,\zeta)$ extends to a 
$C^\infty$ smooth function on $\OD^2\setminus\delta(\T)$, and the norm of
$k_a^\o(\cdot,\zeta)$ in $L^2_a(\o)$ is $\asymp (1-|\zeta|^2)^{-2}$, since
the space is just the usual $L^2_a$ Bergman space with a different equivalent
norm. If $\Lambda$, $\Lambda\subset\D\setminus\{0\}$, is a subsequence of a 
zero sequence (counting multiplicities) for $L^2_a(\o)$, the {\sl extremal 
function} for $\Lambda$ in $L^2_a(\o)$, denoted $G_\Lambda$, is the function 
in $L^2_a(\o)$ that vanishes (counting multiplicities) on $\Lambda$, and has 
biggest positive value at the origin among all such functions of norm $1$.
If $\Lambda$ only contains one point $\lambda$, then we write $G_\lambda$ in
place of $G_\Lambda$. A calculation reveals that
$$G_\lambda(z)=\big(1-\|k_a^\o(\cdot,\lambda)\|^{-2}\big)^{-1/2}
\big(1-\|k_a^\o(\cdot,\lambda)\|^{-2}k_a^\o(z,\lambda)\big),\qquad z\in\D.
\tag3-1$$
Let us agree to say that the one-point zero extremal functions $G_\lambda$
for $L^2_a(\o)$ are asymptotically expansive multipliers if
$$\liminf_{|\lambda|\to1}(1-|\lambda|^2)^{-2}\big(\|G_\lambda f\|_{L^2(\o)}-
\|f\|_{L^2(\o)}\big)\ge0.$$
For instance, if they are expansive for all $\lambda$ in an annulus $r<
|\lambda|<1$, then they are asymptotically expansive.

The extremal functions for subzero sequences (a subzero sequence is a
subsequence of a zero sequence) form a subclass among the 
$L^2_a(\o)$-inner functions $G$, which are described by their property that
$$h(0)=\int_\D h(z)|G(z)|^2\o(z)dS(z),\qquad h\in L^\infty_h(\D).\tag 3-2$$
Note the similarity between \thetag{3-2} and \thetag{1-1}. It is standard to
express these conditions as saying that $\o dS$ and $|G|^2\o dS$ are 
representing measures for the origin.

\proclaim{Theorem 3.1} Suppose that the one-point zero extremal functions in 
$L^2_a(\o)$ are asymptotically expansive multipliers on $L^2_a(\o)$. Then if 
$G$ is $L^2_a(\o)$-inner and continuous on $\OD$, $|G|\ge1$ on $\T$.
\endproclaim

\demo{Proof} It follows from \thetag{3-1} that
$$\big(\|k_a^\o(\cdot,\lambda)\|^2-1\big)\big(|G_\lambda(z)|^2-1\big)=
\|k^\o_a(\cdot,\lambda)\|^{-2}\big|k_a^\o(z,\lambda)\big|^2-2\Re k^\o_a(z,
\lambda)+1.\tag3-3$$
Since both $|G|^2-1$ and $|G_\lambda|^2-1$ annihilate the bounded harmonic
functions with respect to the inner product of $L^2(\o)$, one sees that
$$\int_\D\big(|G(z)|^2-1\big)\,|G_\lambda|^2\o dS=\int_\D
\big(|G(z)|^2-1\big)\big(|G_\lambda|^2-1\big)\o dS=\int_\D|G(z)|^2
\big(|G_\lambda|^2-1\big)\o dS,$$
whence, in view of \thetag{3-3},
$$\multline
\big\|k^\o_a(\cdot,\lambda)\big\|^{-2}\int_\D\big(|G(z)|^2-1\big)\big|
k^\o_a(z,\lambda)\big|^2\o(z)dS(z)\\
=\big(\|k_a^\o(\cdot,\lambda)\|^2-1\big)\int_\D|G(z)|^2\big(|G_\lambda|^2
-1\big)\o dS.
\endmultline$$
As $\lambda$ approaches some point $\lambda_0\in\T$, the left hand side tends
to $|G(\lambda_0)|^2-1$, and the right hand side has a limes inferior that is
nonnegative, by the condition of asymptotic expansion. The assertion of the 
theorem follows.
\qed\enddemo

It has been known for some time that an expansive multiplier that extends
continuously to the closed unit disk necessarily has boundary values of modulus
$\ge1$. The following offers a converse to that statement for one-point zero
set extremal functions. For an analytic function $f$ in $\D$, $Z(f)$ denotes
its zero sequence, counting multiplicities.

\proclaim{Theorem 3.2} Suppose that for some $\lambda\in\D\setminus\{0\}$,
$|G_\lambda(z)|\ge1$ for all $z\in\T$. If also $Z(G_\lambda)=\{\lambda\}$, 
then $G_\lambda$ is an expansive multiplier.
\endproclaim

\demo{Proof} For $\lambda\in\D\setminus\{0\}$, let $L_\lambda$ be the operator
$$L_\lambda f(z)=\frac{f(z)-f(\lambda)}{G_\lambda(z)},\qquad z\in\D\setminus
\{\lambda\},$$
which takes $L^2_a(\o)$ into itself. Note that since $G_\lambda$ is bounded 
and $\o$ is so regular, \thetag{3-2} extends to all $h\in L^2_h$. The functions
$|f|^2$ and $|f-f(\lambda)|^2=|G_\lambda L_\lambda f|^2$ differ by a function
in $L^2_h$, and so
$$\multline
0\le\int_\D\big(|G_\lambda(z)|^2-1\big)^2|L_\lambda f(z)|^2\o(z)dS(z)\\
=\int_\D\big(|G_\lambda(z)|^2-1\big)\big(|f(z)-f(\lambda)|^2-|L_\lambda f(z)|^2
\big)\o(z)dS(z)\\
\int_\D\big(|G_\lambda(z)|^2-1\big)\big(|f(z)|-|L_\lambda f(z)|^2\big)\o(z)
dS(z)\\
=\big(\|G_\lambda f\|^2_{L^2(\o)}-\|f\|^2_{L^2(\o)}\big)-\big(
\|G_\lambda L_\lambda f\|^2_{L^2(\o)}-\|L_\lambda f\|^2_{L^2(\o)}\big).
\endmultline$$
If we apply the above formula to $L_\lambda^j$ in place of $f$, and then sum
over $j=1,\ldots,n-1$, the right hand side becomes a telescoping sum, and the
result is
$$\multline
0\le\int_\D\big(|G_\lambda(z)|^2-1\big)^2\sum_{j=1}^{n-1}|L^j_\lambda f(z)|^2
\o(z)dS(z)\\
=\big(\|G_\lambda f\|^2_{L^2(\o)}-\|f\|^2_{L^2(\o)}\big)-\big(
\|G_\lambda L^{n}_\lambda f\|^2_{L^2(\o)}-\|L^n_\lambda f\|^2_{L^2(\o)}\big).
\endmultline\tag3-4$$
As we let $n$ tend to infinity, the following is obtained.

\proclaim{Proposition 3.3} If $f\in L^2_a(\o)$ is such that $L^n_\lambda f$
tends to $0$ in norm as $n\to+\infty$, then $\|G_\lambda f\|_{L^2(\o)}\ge
\|f\|_{L^2(\o)}$. Indeed,
$$\|G_\lambda f\|^2_{L^2(\o)}=\|f\|_{L^2(\o)}^2+\int_\D
\big(|G_\lambda(z)|^2-1\big)^2\sum_{j=1}^{\infty}|L^j_\lambda f(z)|^2
\o(z)dS(z).$$
\endproclaim

\demo{Remark} If $f=p(G_\lambda)$, where $p$ is a polynomial, then for 
sufficiently large $n$, $L_\lambda^n f=0$. If functions of this type are
dense in $L^2_a(\o)$, then $G_\lambda$ is an expansive multiplier on 
$L^2_a(\o)$. This fact was shown by Serge{\u\i} Shimorin in \cite{\Shi}.
Actually, our result grew out of an effort to generalize Shimorin's method.
\enddemo

We continue the proof of Theorem 3.2. In view of Proposition 3.3, it suffices
to show that the assumption that $G_\lambda$ is greater than or equal to $1$ 
in modulus on $\T$ entails that $L^n_\lambda f$ tends to $0$ in norm as $n\to
+\infty$, for every $f\in L^2_a(\o)$. An approximation argument shows that we
need only do this for polynomials $f$. Let $H^2_\lambda$ be the usual Hardy 
space $H^2$ on the unit disk, endowed with the norm that makes the mapping
$f\mapsto f\circ\varphi_\lambda$ an isometry $H^2\to H^2_\lambda$; here,
$\varphi_\lambda(z)=(\lambda-z)/(1-\bar\lambda z)$ is the standard M\"obius
automorphism. This renormalization has the virtue that $f-f(\lambda)$ is an
orthogonal projection of $f$, and hence has smaller norm than $f$ in 
$H^2_\lambda$. Division by $G_\lambda$ then reduces the norm even more, and
hence $\|L_\lambda f\|_{H^2_\lambda}\le \|f\|_{H^2_\lambda}$. It follows that
$\|L^n_\lambda f\|_{H^2_\lambda}\le \|f\|_{H^2_\lambda}$ for all $n=1,2,3,
\ldots$. The functions $L^n_\lambda f$ clearly form a normal family on $\D$,
so that we can select a normal limit $g\in H^2_\lambda$. This normal limit
has $\|L_\lambda g\|_{H^2_\lambda}=\|g\|_{H^2_\lambda}$, which is only 
possible if both the projection $g\mapsto g-g(\lambda)$ and the division of
$g-g(\lambda)$ by $G_\lambda$ are isometric. Unless $g$ is the $0$ function,
the division is isometric only if $|G_\lambda(z)|=1$ everywhere on $\T$. By 
the maximum principle, and the fact that $G_\lambda$ is nonconstant, we obtain
$|G_\lambda(z)|<1$ throughout $\D$. This, however, makes \thetag{3-2} 
impossible for $G=G_\lambda$ and $h=1$. The only remaining logical alternative
is that $g=0$. But then every normal limit of $L^n_\lambda f$ is $0$, whence
$L^n_\lambda f$ tends to $0$ as $n\to+\infty$, uniformly on compact subsets
of $\D$. Since the functions $L^n_\lambda f$ were known to be bounded in
$H^2_\lambda$, an elementary argument now shows that $L^n_\lambda f\to0$ in 
the norm of $L^2_a$, which is equivalent to that of $L^2_a(\o)$. The proof
is complete.
\qed\enddemo

\heading 4. A refinement of Shimorin's theorem on compositions
\endheading

The following result is a slight sharpening of a theorem of Shimorin
\cite{\bf\Shi}. The space $L^p_a(\C,\mu)$ denotes the closure of the
polynomials in the space $L^p(\C,\mu)$, and $\|\cdot\|_{L^p(\mu)}$ stands
for the norm in the latter space. Similarly, $L^\infty_a(\C,\mu)$ is the
weak-star closure of the polynomials in $L^\infty(\C,\mu)$.

\proclaim{Theorem 4.1} Let $p$ be a positive even integer, and let $\mu$ be
a compactly supported Borel probability measure in $\C$ such that
$$q(0)=\int_\C q(z)\,d\mu(z)$$
holds for all polynomials $q$. Suppose $\varphi\in L^\infty_a(\C,\mu)$
satisfies
$$q(0)=\int_\C q(z)\,|\varphi(z)|^p\,d\mu(z)$$
for all polynomials $q$. Then
$$\|q\circ\varphi\|_{L^p(\mu)}\le\|(q\circ\varphi)\varphi\|_{L^p(\mu)},$$
again for all polynomials $q$.
\endproclaim

\demo{Proof} Let $\calP_a$ stand for the linear space of all polynomials,
and let ${\calP}_h$ stand for the space of harmonic polynomials (that is,
all functions of type $q_1(z)+\overline{q}_2(z)$, with $q_1,q_2\in{\calP}_a$).
It is clear from the hypotheses and the standard approximation of bounded
harmonic functions in $\D$ by harmonic polynomials that
$$\int_\C u(z)\big(|\varphi(z)|^p-1\big)\,d\mu(z)=0,\qquad u\in L^\infty_h
(\D).\tag4-1$$
For $f\in{\calP}_a$, we define
$$Tf(z)=\frac{f(z)-f(0)}{z},$$
which is another polynomial; the mapping $T$ is known as the backward shift.
Since $zTf(z)=f(z)-f(0)$, we see that
$$\Delta|zTf(z)|^2=\Delta|f(z)|^2.\tag4-2$$
We write $p=2N$, where $N$ is a positive integer, and let $q\in{\calP}_a$.
Also, let $g(z)=q(z)^N$, and observe that as $g$ is a polynomial, the sum
$$\sum_{n=1}^\infty|T^n g(z)|^2$$
is actually finite. Repeated applications of \thetag{4-2} show that
$$\Delta\left((1-|z|^2)\sum_{n=1}^\infty|T^n g(z)|^2\right)=-\Delta|g(z)|^2,$$
whence
$$u_g(z)=|g(z)|^2+(1-|z|^2)\sum_{n=1}^\infty|T^n g(z)|^2$$
is harmonic; in fact, it is in ${\calP}_h$. We now see that, in view of
\thetag{4-1},
$$\multline
\|(q\circ\varphi)\varphi\|_{L^p(\mu)}-\|q\circ\varphi\|_{L^p(\mu)}=
\int_\C \big|q(\varphi(z))\big|^p\big(|\varphi(z)|^p-1\big)\,d\mu(z)\\
=\int_\C \big|g(\varphi(z))\big|^2\big(|\varphi(z)|^{2N}-1\big)\,d\mu(z)
=\int_\C \big|g(\varphi(z))\big|^2\big(|\varphi(z)|^{2N}-1\big)\,d\mu(z)\\
=\int_\C \left(\big|g(\varphi(z))\big|^2-u_g(\varphi(z))\right)
\big(|\varphi(z)|^{2N}-1\big)\,d\mu(z)\\
=\int_\C \sum_{n=1}^\infty\big|(T^n g)(\varphi(z))\big|^2
\big(|\varphi(z)|^{2}-1\big)\big(|\varphi(z)|^{2N}-1\big)\,d\mu(z)\ge0,
\endmultline$$
where the last relation holds simply because the integrand is greater than
or equal to $0$.
The proof is complete.
\qed\enddemo

We now show that the assumption in Theorem 4.1 that $p$ is an even integer
cannot be avoided.

\proclaim{Theorem 4.2} Let $p$ be a positive real number, other than an
even integer. Then there is a Borel probability measure supported in the
closed unit disk $\OD$, with the following properties:
\roster
\item{}$\quad u(0)=\int_{\OD}u\,d\mu$ for all $u\in{\calP}_h$,
\item{}$\quad u(0)=\int_{\OD}u(z)\,|2z|^p\,d\mu$ for all $u\in{\calP}_h$,
\item{}$\quad \|(z-\frac12)(2z)\|_{L^p(\mu)}<\|(z-\frac12)\|_{L^p(\mu)}$.
\endroster
\endproclaim

\demo{Proof} Consider the real Banach space $C(\OD)$ consisting of the
real-valued continuous functions on $\OD$, with the supremum norm. The
dual space of $C(\OD)$ is identified with the space of real-valued Borel
measures on $\OD$, and we will show the existence of the measure $\mu$ by
showing that a suitable continuous linear functional on $C(\OD)$ exists.
Let $L$ be the linear subspace of $C(\OD)$ consisting of the functions of
the form
$$u(z)+\big(|2z|^p-1)\,v(z),\qquad u,v\in{\calP}_h,\tag{4-3}$$
and let $\Lambda:L\to\R$ be the linear functional on $L$ taking the
functions of the form \thetag{4-3} to $u(0)$.

Since
$$\multline
|u(0)|\le\sup\big\{|u(z)|:\,|z|=\tfrac12\big\}=
\sup\big\{\big|u(z)+(|2z|^p-1)v(z)\big|:\,|z|=\tfrac12\big\}\\
\le\sup\big\{\big|u(z)+(|2z|^p-1)v(z)\big|:\,|z|\le1\big\},
\endmultline$$
and $\Lambda(1)=1$, we see that the norm of $\Lambda$ as a linear functional
from the subspace $L$ of $C(\OD)$ to $\R$ equals $1$.
We wish to extend $\Lambda$ in a norm-preserving way to $L_1$, the subspace
of $C(\OD)$ spanned by $L$ and the function
$$F(z)=\big(|2z|^p-1\big)\,|z-\tfrac12|^p.$$
Suppose $\Lambda_1:L_1\to\R$ is defined by
$$\Lambda_1(f+\lambda F)=\Lambda(f)+\lambda\alpha,\qquad f\in L,\,\,
\lambda\in\R,$$
where $\alpha\in\R$. The functional $\Lambda_1$ is obviously a linear
extension of $\Lambda$, and by the proof of the Hahn-Banach theorem (see
\cite{\Rud}), the norm of $\Lambda_1$ will be $1$ provided that
$$\sup\big\{\Lambda(f)-\|f-F\|_{C(\OD)}:\,f\in L\big\}\le\alpha\le
\inf\big\{\Lambda(f)+\|f-F\|_{C(\OD)}:\,f\in L\big\};$$
it is part of the Hahn-Banach theorem that the left hand side is less than or
equal to the right hand side, so that a permissible $\alpha$ can be found.
We shall choose $\alpha$ to be the quantity on the left, that is,
$$\alpha=\sup\Big\{u(0)-\sup\big\{\big|u(z)+(|2z|^p-1)v(z)-(|2z|^p-1)
|z-\tfrac12|^p\big|:\,|z|\le1\big\}:\,\,u,v\in{\calP}_h\Big\}.\tag{4-4}$$
Again by the Hahn-Banach theorem, there is a norm-preserving extension of
$\Lambda_1$ to a bounded linear functional on $C(\OD)$. Identifying this
functional with a real-valued Borel measure $\mu$ on $\OD$, we see that it
has the properties
$$\gather
u(0)=\int_{\OD}\big(u(z)+(|2z|^p-1)v(z)\big)\,d\mu(z),
\qquad u,v\in{\calP}_h,\tag{4-5}\\
\int_{\OD}\big(|2z|^p-1\big)|z-\tfrac12|^p\,d\mu(z)=\alpha,\tag{4-6}\\
\|\mu\|=1.\tag{4-7}
\endgather$$
By \thetag{4-5}, $\mu$ has properties (1) and (2). As we combine \thetag{4-5}
with \thetag{4-7}, we see that
$$\Lambda(1)=\int_{\OD}d\mu=\int_{\OD}d|\mu|=\|\mu\|=1,$$
so that $\mu$ is a probability measure. We shall complete the proof by
demonstrating that $\alpha<0$, thus showing that property (3) is satisfied.

Suppose on the contrary that $\alpha\ge0$. From \thetag{4-4} we see that
there must exist $u_n,v_n\in{\calP}_h$ for $n=1,2,3,\ldots$, such that
$$\sup_{z\in\OD}\Big|u_n(z)+\big(|2z|^p-1\big)\,v_n(z)-
\big(|2z|^p-1\big)\,|z-\tfrac12|^p\Big|<u_n(0)+\frac1n.$$
It follows from this that
$$u_n(z)-u_n(0)+\big(|2z|^p-1\big)\,v_n(z)-
\big(|2z|^p-1\big)\,|z-\tfrac12|^p<\frac1n,\qquad z\in\OD.\tag{4-8}$$
Plugging in $z=0$ into this gives
$$v_n(0)>2^{-p}-\frac1n.$$
We now set
$$\gather
f_n(z)=u_n(z)-u_n(0)+(2^p-1)v_n(z),\\
g_n(z)=u_n(z)-u_n(0)+((\tfrac32)^p-1\big)v_n(z).
\endgather$$
From \thetag{4-8} we can see that
$$u_n(z)-u_n(0)+\big(|2z|^p-1\big)\,v_n(z)<\big(|2z|^p-1\big)(|z|+\tfrac12)^p
+\frac1n,\qquad\tfrac12\le|z|\le1.$$
Plugging $|z|=1$ and $|z|=\frac34$ into this yields
$$f_n(z)<(2^p-1)(\tfrac32)^p+\frac1n,\qquad |z|=1,\tag{4-10}$$
and
$$g_n(z)<((\tfrac32)^p-1)(\tfrac54)^p+\frac1n,\qquad |z|=\tfrac34.
\tag{4-11}$$
By \thetag{4-9}, we see that
$$\gather
f_n(0)>(2^p-1)(2^{-p}-n^{-1}),\tag{4-12}\\
g_n(0)>\big((\tfrac32)^p-1\big)(2^{-p}-n^{-1}).\tag{4-13}
\endgather$$
The functions $f_n$ and $g_n$ are harmonic polynomials. It follows from
\thetag{4-10}--\thetag{4-13} that $\{f_n\}_n$ and $\{g_n\}_n$ form normal
families of harmonic functions in $|z|<\frac34$, and hence there are
subsequences $\{f_{n_j}\}_j$ and $\{g_{n_j}\}_j$ that converge uniformly on
compact subsets of $|z|<\frac34$. It follows that the corresponding
subsequences of $\{u_n\}_n$ and $\{v_n\}_n$, which may be represented by
$$\gather
u_{n_j}(z)-u_{n_j}(0)=\big(2^p-(\tfrac32)^p\big)^{-1}
\big( (2^p-1)g_{n_j}(z)-((\tfrac32)^p-1)f_{n_j}(z) \big),\\
v_{n_j}(z)=\big(2^p-(\tfrac32)^p\big)^{-1}\big(f_{n_j}(z)-g_{n_j}(z)\big),
\endgather$$
converge uniformly on compact subsets of $|z|<\frac34$. Say $u_{n_j}(z)
-u_{n_j}(0)\to U(z)$ and $u_{n_j}(z)\to V(z)$ as $j\to+\infty$, where $U(z)$
and $V(z)$ are harmonic in $|z|<\frac34$. It follows
from \thetag{4-8} that
$$U(z)+\big(|2z|^p-1\big)V(z)-\big(|2z|^p-1\big)|z-\tfrac12|^p\le0,\qquad
|z|<\tfrac34.\tag{4-14}$$
Hence $U(z)\le0$ for $|z|=\frac12$. Since $U$ is harmonic in $|z|<\frac34$
and $U(0)=0$, we conclude that $U(z)\equiv0$. The relation \thetag{4-14} now
reads
$$\big(|2z|^p-1\big)\big(V(z)-|z-\tfrac12|^p\big)\le0,\qquad
|z|<\tfrac34.\tag{4-15}$$
However, for $|z|<\frac12$, \thetag{4-15} says that
$$V(z)-|z-\tfrac12|^p\ge0,$$
and for $|z|>\frac12$, it says that
$$V(z)-|z-\tfrac12|^p\le0.$$
It follows by continuity that
$$V(z)=|z-\tfrac12|^p,\qquad |z|=\tfrac12.$$
This is a contradiction, since the function $|z-\frac12|^p$ fails to be
infinitely differentiable at $z=\frac12$.
\qed\enddemo

%\endRefs


\Refs

\ref\no\Hed
\by H.~Hedenmalm
\paper A factorization theorem for square area-integrable analytic functions
\jour J. Reine Angew. Math.\vol422\yr 1991\pages 45--68
\endref

\ref\no\Rud
\by W.~Rudin
\book Real and Complex Analysis, 3rd Edition
\publ McGraw-Hill Book Co\yr1987
\endref

\ref\no\Shi
\by S.~Shimorin
\paper Single point extremal functions in weighted Bergman spaces.
Nonlinear boundary-value problems and some questions of function theory
\jour J. Math. Sci. \vol80{\rm, no. 6} \yr1996\pages 2349--2356
\endref
\enddocument